\documentclass{amsart}
\theoremstyle{plain}

\newtheorem{theorem}{Theorem}

\newcommand{\Z}{\mathbb{Z}}
\newcommand{\N}{\mathbb{N}}
\DeclareMathOperator{\tr}{tr}
\DeclareMathOperator{\Ima}{Im}
\newcommand{\Ca}{{\mathbf{C}}}

\begin{document}
\title[On matrices in prescribed conjugacy classes]{On matrices in
prescribed conjugacy classes with no common invariant subspace and sum zero}
\author{William Crawley-Boevey}
\address{Department of Pure Mathematics, University of Leeds, Leeds LS2 9JT, UK}
\email{w.crawley-boevey@leeds.ac.uk}
\thanks{Mathematics Subject Classification (2000): Primary 15A24, 16G20; Secondary 34M50.}

\begin{abstract}
We determine those $k$-tuples of conjugacy classes of matrices, from which it is possible
to choose matrices which have no common invariant subspace and have sum zero.
This is an additive version of the Deligne-Simpson problem.
We deduce the result from earlier work of ours on preprojective algebras
and the moment map for representations of quivers.
Our answer depends on the root system for a Kac-Moody Lie algebra.
\end{abstract}
\maketitle

\section{Introduction}
A problem considered by Deligne, Simpson \cite{Simp}, and others,
concerns the equation
\begin{equation}
\label{eq1}
A_1 A_2 \dots A_k = 1
\end{equation}
for $n\times n$ matrices $A_i$, with entries in an algebraically closed field $K$, say.
The additive analogue is the equation
\begin{equation}
\label{eq2}
A_1+A_2 + \dots+A_k = 0.
\end{equation}
A solution of the additive equation over the field of
complex numbers leads to a solution of the multiplicative
equation, as the monodromy of a Fuchsian system of ODEs.
See for example \cite{AB}.

The problem is to describe those $k$-tuples of conjugacy classes $C_1,\dots,C_k$
for which there exists a solution to (\ref{eq1}) or (\ref{eq2}) with $A_i\in C_i$,
and which is \emph{irreducible}, meaning that there is no nonzero proper
subspace of $K^n$ which is invariant for all the $A_i$.
It is stated in this form, and studied, by Kostov
\cite{Kostovmon,KostovCR,Kostovaspects,Kostovnext},
who calls it the `Deligne-Simpson problem'.

In this paper we solve the problem for equation (\ref{eq2}).
It is a consequence of our work on preprojective algebras and the
moment map for representations of quivers~\cite{CBce,CBmm,CBH}.

For each $1\le i\le k$, choose elements $\xi_{i,1},\xi_{i,2},\dots,\xi_{i,d_i}\in K$
(not necessarily distinct) with
\begin{equation}
\label{e:prodxis}
\prod_{j=1}^{d_i} (A - \xi_{i,j}1) = 0
\end{equation}
for the matrices $A\in C_i$.
For notational convenience we assume that $d_i\ge 2$ for all $i$.
The conjugacy class $C_i$ is uniquely
determined by the elements $\xi_{i,j}$
and the ranks $r_{i,j}$ of the partial products $\prod_{\ell=1}^j (A - \xi_{i,\ell}1)$.
Clearly
\[
n = r_{i,0} \ge r_{i,1} \ge \dots \ge r_{i,d_i} = 0
\]
and
\begin{equation}
\label{e:rcond}
r_{i,j-1} - r_{i,j} \ge r_{i,\ell-1} - r_{i,\ell}
\end{equation}
whenever $1\le j<\ell\le d_i$ and $\xi_{i,j} = \xi_{i,\ell}$.

Our answer depends on the root system for the Kac-Moody Lie
algebra with symmetric generalized Cartan matrix $\Ca$ whose diagram is

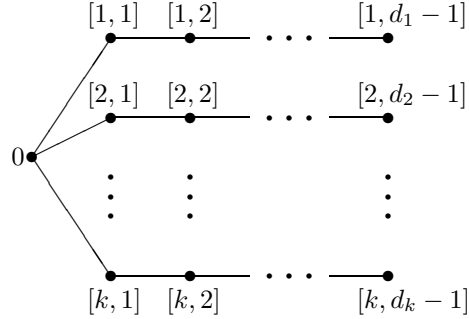
\setlength{\unitlength}{1.5pt}
\[
\begin{picture}(110,80)
\put(10,40){\circle*{2.5}}
\put(30,10){\circle*{2.5}}
\put(30,50){\circle*{2.5}}
\put(30,70){\circle*{2.5}}
\put(50,10){\circle*{2.5}}
\put(50,50){\circle*{2.5}}
\put(50,70){\circle*{2.5}}
\put(100,10){\circle*{2.5}}
\put(100,50){\circle*{2.5}}
\put(100,70){\circle*{2.5}}
\put(10,40){\line(2,-3){20}}
\put(10,40){\line(2,1){20}}
\put(10,40){\line(2,3){20}}
\put(30,10){\line(1,0){35}}
\put(30,50){\line(1,0){35}}
\put(30,70){\line(1,0){35}}
\put(85,10){\line(1,0){15}}
\put(85,50){\line(1,0){15}}
\put(85,70){\line(1,0){15}}
\put(70,10){\circle*{1}}
\put(75,10){\circle*{1}}
\put(80,10){\circle*{1}}
\put(70,50){\circle*{1}}
\put(75,50){\circle*{1}}
\put(80,50){\circle*{1}}
\put(70,70){\circle*{1}}
\put(75,70){\circle*{1}}
\put(80,70){\circle*{1}}
\put(30,25){\circle*{1}}
\put(30,30){\circle*{1}}
\put(30,35){\circle*{1}}
\put(50,25){\circle*{1}}
\put(50,30){\circle*{1}}
\put(50,35){\circle*{1}}
\put(100,25){\circle*{1}}
\put(100,30){\circle*{1}}
\put(100,35){\circle*{1}}
\put(5,38){0}
\put(24,2){$[k,1]$}
\put(24,54){$[2,1]$}
\put(24,74){$[1,1]$}
\put(44,2){$[k,2]$}
\put(44,54){$[2,2]$}
\put(44,74){$[1,2]$}
\put(92,2){$[k,d_k-1]$}
\put(92,54){$[2,d_2-1]$}
\put(92,74){$[1,d_1-1]$}
\end{picture}
\]

\noindent
Thus the rows and columns of $\Ca$ are indexed by the vertex set
\[
I = \{ 0 \} \cup \{ [i,j] \mid 1\le i\le k, 1\le j \le d_i-1\},
\]
the diagonal entries are $\Ca_{vv}=2$ for $v\in I$, and the off-diagonal entries are $-1$ if there is an edge joining the
two vertices, otherwise zero.

The root system can be considered as a subset of the set of column vectors $\Z^I$.
It includes the \emph{coordinate vectors} $\epsilon_v$
($v\in I$) and the \emph{fundamental region}, which consists of the nonzero elements $\alpha\in\N^I$ which have
connected support and all components of $\Ca \alpha$ $\le 0$.
It is then closed up under change of
sign and the action of the Weyl group, which is generated by the reflections $s_v:\Z^I\to \Z^I$,
defined for $v\in I$ by $s_v(\alpha) = \alpha - (\epsilon_v^T \Ca \alpha)\epsilon_v$.
The roots coming from a coordinate vector are called \emph{real roots}; those coming from an element of
the fundamental region are called \emph{imaginary roots}.

If $\lambda\in K^I$ we denote by $R^+_\lambda$ the set of positive
roots $\alpha\in\N^I$ with the property that $\lambda\cdot\alpha := \sum_{v\in I} \lambda_v \alpha_v = 0$.
For $\alpha\in\Z^I$ we define $p(\alpha) = 1 - \frac12 \alpha^T \Ca \alpha \in \Z$. For a root $\alpha$,
one has $p(\alpha)\ge 0$, with equality if and only if $\alpha$ is a real root.
We denote by $\Sigma_\lambda$ the set of $\alpha\in R^+_\lambda$
with the property that
\begin{equation}
\label{e:pmax}
p(\alpha) > p(\beta^{(1)}) + p(\beta^{(2)}) + \dots
\end{equation}
for any decomposition $\alpha = \beta^{(1)} + \beta^{(2)}+ \dots$ as a sum of
two or more elements of $R^+_\lambda$.
Other characterizations of this set are given by
Theorems 5.6 and 8.1 of \cite{CBmm}.
Our result is as follows.

\begin{theorem}
\label{t:main}
There is an irreducible solution to equation {\rm (\ref{eq2})} with matrices $A_i\in C_i$
if and only if $\alpha\in\Sigma_\lambda$,
where $\alpha$ is defined by $\alpha_0 = n$ and $\alpha_{[i,j]} = r_{i,j}$,
and $\lambda$ is defined by $\lambda_0 = - \sum_{i=1}^k \xi_{i,1}$ and
$\lambda_{[i,j]} = \xi_{i,j} - \xi_{i,j+1}$. In this case, if $\alpha$ is a real root,
then any other solution to {\rm (\ref{eq2})} with matrices in $C_i$ is conjugate
to this solution, while if $\alpha$ is an imaginary root, then there are infinitely
many non-conjugate irreducible solutions.
\end{theorem}

We prove Theorem \ref{t:main} in Section \ref{s:proof}.
The case where there is an irreducible solution which is unique up
to conjugacy is called the \emph{rigid} case, and we discuss it in
Section~\ref{s:rigid}.
Another special case worth considering is when the $C_i$ are nilpotent
conjugacy classes. We discuss it in Section~\ref{s:nilpotent}.
Finally, at the opposite extreme, we discuss the case when
the $C_i$ have generic eigenvalues in Section~\ref{s:generic}.
I would like to thank L. Hille for bringing this problem to my attention,
and several referees for their helpful remarks.

\section{Deformed preprojective algebras}
In this section we recall the deformed
preprojective algebras of \cite{CBH}, and a result
proved in \cite{CBmm}.

Let $Q$ be a quiver with vertex set $I$.
Let $\Ca$ be the corresponding generalized Cartan matrix,
with rows and columns indexed by $I$, so with diagonal entries $\Ca_{vv}=2$
and off-diagonal entries $\Ca_{vw} = -(n_{vw}+n_{wv})$
where $n_{vw}$ denotes the number of arrows with head at $v$ and tail at $w$.
Let $R_\lambda^+$, $p(\alpha)$ and $\Sigma_\lambda$ be defined
as in the introduction.
(For simplicity we are assuming that
$Q$ has no loops at vertices, otherwise these definitions need to be modified.)

The \emph{double} of $Q$ is the quiver $\overline{Q}$ obtained from $Q$ by adjoining
a reverse arrow $a^*:w\to v$ for each arrow $a:v\to w$ in $Q$.
If $\lambda\in K^I$, the \emph{deformed preprojective algebra} is
\[
\Pi^\lambda = K\overline{Q} / \bigl( \sum_{a\in Q} (aa^*-a^*a) - \sum_{v\in I} \lambda_v e_v \bigr)
\]
where $K\overline{Q}$ is the path algebra of $\overline{Q}$, and $e_v$ denotes the
trivial path at vertex $v$.

Recall that a \emph{representation} $X$ of $\overline{Q}$ is given by a vector space
$X_v$ for each vertex and a linear map $a:X_v\to X_w$ for each arrow $a:v\to w$ in $\overline{Q}$.
The \emph{dimension vector} of $X$ is the vector in $\N^I$ whose components are the
dimensions of the spaces $X_v$.

Representations of $\Pi^\lambda$ are the same as representations of $\overline{Q}$
in which the linear maps satisfy
\[
\sum_{\substack{a\in Q \\ h(a)=v}} aa^* - \sum_{\substack{a\in Q \\ t(a)=v}} a^*a = \lambda_v 1,
\]
for each $i$.
Here $h(a)$ and $t(a)$ denote the head and tail vertices of an arrow $a$.
The concepts of `isomorphism' of representations and `simple' representations
are straightforward.

In \cite[Theorem 1.2]{CBmm} we have proved the following.

\begin{theorem}
\label{t:quoted}
There is a simple representation of $\Pi^\lambda$ of dimension vector $\alpha$ if and only
if $\alpha\in\Sigma_\lambda$.
If $\alpha$ is a real root, the simple representation is unique
up to isomorphism, and is the only representation of $\Pi^\lambda$ of dimension vector $\alpha$.
If $\alpha$ is an imaginary root there are infinitely many non-isomorphic simple representations.
\end{theorem}

Some parts are not explicitly stated there, but they follow easily:
the number of simple representations is discussed in the remarks following
the statement of the theorem; and if there is a non-simple representation $X$ whose
dimension vector is a real root $\alpha$, then the decomposition of $\alpha$ as the sum
of the dimension vectors of the composition factors of $X$ shows that $\alpha\notin\Sigma_\lambda$.

\section{Proof of Theorem \ref{t:main}}
\label{s:proof}
Let $\Ca$ be the generalized Cartan matrix constructed in the introduction,
let  $Q$ be the quiver obtained from its diagram
by orienting all arrows towards the vertex $0$, and denote by $a_{i,j}$ the arrow with tail at
vertex $[i,j]$ and head at $[i,j-1]$.
(For convenience we use the convention that $[i,0]$ denotes the vertex $0$ for any $i$.)
Let $\alpha$ and $\lambda$ be defined as in the statement of Theorem \ref{t:main}.

Given a solution to equation (\ref{eq2}) with matrices $A_i\in C_i$,
we construct a representation $X$ of $\Pi^\lambda$
of dimension vector $\alpha$ in which all the linear maps $a_{[i,j]}$ are injective and $a_{[i,j]}^*$
are surjective. The vector spaces at each vertex are
$X_0 = K^n$ and $X_{[i,j]} = \Ima (\prod_{\ell=1}^j (A_i - \xi_{i,\ell}1))$.
The linear map for arrow $a_{i,j}$ is the inclusion, and the linear map
for $a_{i,j}^*$ is
\[
a_{i,j}^* = (A_i - \xi_{i,j} 1)|_{X_{[i,j-1]}}.
\]
This is well-defined since $(A_i - \xi_{i,j} 1)(X_{[i,j-1]}) = X_{[i,j]}$.

We show that $X$ is a representation of $\Pi^\lambda$. For $j < d_i-1$ we have
\[
a_{i,j+1} a_{i,j+1}^* - a_{i,j}^* a_{i,j} = (A_i - \xi_{i,j+1}1) - (A_i - \xi_{i,j}1) = \lambda_{[i,j]} 1
\]
which is the relation at vertex $[i,j]$. Also,
since the restriction of $A_i - \xi_{i,d_i}1$ to $X_{[i,d_i-1]}$ is zero
by (\ref{e:prodxis}), we have
\[
- a_{i,d_i-1}^* a_{i,d_i-1} = - (A_i - \xi_{i,d_i-1} 1) = \lambda_{[i,d_i-1]} 1
\]
which is the relation at $[i,d_i-1]$. Finally (\ref{eq2}) implies that
\[
\sum_{i=1}^k a_{i,1} a_{i,1}^* = \sum_{i=1}^k (A_i - \xi_{i,1} 1) = \lambda_0 1
\]
which is the relation at 0.

Conversely, given a representation $X$ of $\Pi^\lambda$ of dimension vector $\alpha$
in which all the linear maps $a_{[i,j]}$ are injective and $a_{[i,j]}^*$
are surjective, we construct a solution to equation (\ref{eq2}) with matrices $A_i\in C_i$.
Identify $X_0$ with $K^n$.
Define $A_i = a_{i,1} a^*_{i,1} + \xi_{i,1} 1$.
For $j<d_i-1$ the relation at vertex $[i,j]$ is
$a_{i,j+1} a_{i,j+1}^* - a_{i,j}^* a_{i,j} = \lambda_{[i,j]} 1$,
which implies that
\begin{equation}
\label{e:rone}
a^*_{i,j} (a_{i,j} a^*_{i,j} +(\xi_{i,j}+c) 1) = (a_{i,j+1} a^*_{i,j+1} + (\xi_{i,j+1}+c) 1) a^*_{i,j}
\end{equation}
for $c\in K$.
The relation at $[i,d_i-1]$ is
$- a^*_{i,d_i-1} a_{i,d_i-1} = \lambda_{[i,d_i-1]} 1$,
which implies that
\begin{equation}
\label{e:rtwo}
a^*_{i,d_i-1} (a_{i,d_i-1} a^*_{i,d_i-1} +(\xi_{i,d_i-1}+c) 1) = (\xi_{i,d_i}+c) a^*_{i,d_i-1}.
\end{equation}
For $j<d_i-1$, equation (\ref{e:rone}) gives
\begin{equation}
\label{e:prodone}
\begin{split}
a^*_{i,j} \dots a^*_{i,2} a^*_{i,1} (A_i + c 1)
&= a^*_{i,j} \dots a^*_{i,2} a^*_{i,1} ( a_{i,1} a^*_{i,1} + (\xi_{i,1} + c) 1)
\\
&= a^*_{i,j} \dots a^*_{i,2} ( a_{i,2} a^*_{i,2} + (\xi_{i,2} + c) 1) a^*_{i,1}
\\
&=\dots
\\
&= a^*_{i,j} (a_{i,j} a^*_{i,j} + (\xi_{i,j} + c) 1) a^*_{i,j-1} \dots a^*_{i,1}
\\
&= (a_{i,j+1} a^*_{i,j+1} + (\xi_{i,j+1}+c)1) a^*_{i,j} \dots a^*_{i,1}.
\end{split}
\end{equation}
Now equation (\ref{e:rtwo}) gives
\begin{equation}
\label{e:prodtwo}
a^*_{i,d_i-1} a^*_{i,d_i-2} \dots a^*_{i,1} (A_i + c1)
= (\xi_{i,d_i}+c)  a^*_{i,d_i-1} a^*_{i,d_i-2} \dots a^*_{i,1}.
\end{equation}
Taking $c = -\xi_{i,j+1}$ in (\ref{e:prodone}) we have
\begin{equation}
\label{e:prodonesp}
a^*_{i,j} \dots a^*_{i,2} a^*_{i,1} (A_i - \xi_{i,j+1} 1)
= a_{i,j+1} a^*_{i,j+1} a^*_{i,j} \dots a^*_{i,1}.
\end{equation}
and hence by induction
\[
\prod_{\ell=1}^j (A_i - \xi_{i,\ell}1) = a_{i,1} a_{i,2} \dots a_{i,j} a^*_{i,j} \dots a^*_{i,2} a^*_{i,1}
\]
for $j\le d_i-1$. Now (\ref{e:prodtwo}) gives
\[
\prod_{\ell=1}^{d_i} (A_i - \xi_{i,\ell}1) = 0.
\]
Using the injectivity of the $a_{i,j}$ and the surjectivity of the $a^*_{i,j}$ these two
product formulas imply that $A_i\in C_i$.
Finally, the relation at vertex 0 is $- \sum_{i=1}^k a_{i,1} a^*_{i,1} = \lambda_0 1$,
which implies that the matrices $(A_i)$ are a solution to equation (\ref{eq2}).

Thus we have a correspondence between solutions of equation
(\ref{eq2}) with matrices $A_i\in C_i$ and representations of $\Pi^\lambda$
of dimension vector $\alpha$ in which all the linear maps $a_{[i,j]}$ are injective and $a_{[i,j]}^*$
are surjective. It is clear that this gives a 1-1 correspondence between conjugacy classes
of solutions and isomorphism classes of representations.
Thus, provided we can show that irreducible solutions correspond to simple representations
of $\Pi^\lambda$, the result follows from Theorem \ref{t:quoted}.

If the solution to (\ref{eq2}) is irreducible, then $X$ as constructed above
is a simple representation of $\Pi^\lambda$,
for if $Y$ is a subrepresentation, the irreducibility implies that $Y_0 = 0$
or $Y_0 = X_0$. But if $Y_0=0$ then $Y=0$ since the linear maps $a_{i,j}$ are all injective,
and if $Y_0 = X_0$ then $Y=X$ since
the linear maps $a^*_{i,j}$ are all surjective.

Conversely, suppose that $X$ is a simple representation of $\Pi^\lambda$
of dimension vector~$\alpha$.
We show that the linear maps $a_{i,j}$ are injective.
For a contradiction, let $x\in X_{[i,\ell]}$ be a nonzero element in the
kernel of $a_{i,\ell}$.
We define elements $x_j\in X_{[i,j]}$ for $j\ge \ell$
by setting $x_\ell = x$ and $x_{j+1} = a^*_{i,j+1}(x_j)$ for $j\ge\ell$.
An induction, using the relation
\[
a_{i,j+1} a_{i,j+1}^* - a_{i,j}^* a_{i,j} = \lambda_{[i,j]} 1
\]
for $j<d_i-1$, shows that $a_{i,j+1}(x_{j+1})$ is a multiple of $x_j$ for $j\ge\ell$.
It follows that the elements $x_j$ span a subrepresentation of $X$. It must be a
proper subrepresentation since it is zero at the vertex 0, but
$\alpha_0 = n \ge r_{i,\ell} \neq 0$.
This contradicts the simplicity of $X$.

A dual argument shows that the linear maps $a^*_{i,j}$ are surjective.

Now suppose that $Y_0$ is an subspace of $X_0 = K^n$ which is invariant under the $A_i$.
Define subspaces $Y_{[i,j]} \subseteq X_{[i,j]}$ via $Y_{[i,j]} = a^*_{i,j}(Y_{[i,j-1]})$
for $j\ge 1$. Thanks to equation (\ref{e:prodonesp})
this defines a subrepresentation of $X$.
But $X$ is simple.
Thus the solution $(A_i)$ is irreducible.

\section{The rigid case}
\label{s:rigid}
An irreducible solution is said to be \emph{rigid} if any other solution (with the
matrices in the same conjugacy classes) is conjugate to it. Two additional problems
which may be considered as part of the `Deligne-Simpson problem' are
to enumerate all cases in which there is a rigid solution,
and to construct the corresponding solutions.

For the multiplicative equation (\ref{eq1}), Katz \cite{Katz} has an algorithm
which enables one to determine whether or not a given $k$-tuple of conjugacy classes has
a rigid solution. Katz does not, however, attempt to enumerate all rigid cases, instead
remarking that `Even a cursory glance \dots leaves one
with the impression that there is a fascinating bestiary waiting to be compiled'.

Katz's algorithm was simplified by Dettweiler and Reiter \cite{DR}, who also adapted
it to the additive equation (\ref{eq2}), but again made no attempt to enumerate
the cases.

Our results show that for equation (\ref{eq2}), there is a rigid solution if and only if
\begin{itemize}
\item[(i)]
$\alpha$ is a positive real root,
\item[(ii)]
$\lambda\cdot\alpha = 0$, and
\item[(iii)]
for any decomposition $\alpha = \beta^{(1)} + \beta^{(2)}+ \dots$
as a sum of two or more positive roots,
one has $\lambda\cdot\beta^{(s)}\neq 0$ for some $s$.
\end{itemize}
Thus we see that the problem of enumerating all rigid
cases amounts to enumerating the positive real roots
$\alpha$, and for each one, determining the set of $\lambda$ satisfying (ii) and (iii).
Note that since a positive real root cannot be proportional to any other positive
root, this set contains the general element of $\{ \lambda\in K^I \mid \lambda\cdot\alpha=0\}$.

For example, given 3 conjugacy classes $C_1$, $C_2$, $C_3$ of $2\times 2$ matrices, there is
a rigid solution to equation (\ref{eq2}) if and only if
\begin{itemize}
\item[(i${}'$)]
none of the $C_i$ consists of a multiple of the identity matrix,
\item[(ii${}'$)]
the sum of all six eigenvalues for the three conjugacy classes is zero, and
\item[(iii${}'$)]
the sum of three eigenvalues, one for each $C_i$, is always nonzero.
\end{itemize}
Namely, taking $\xi_{i,1},\xi_{i,2}$ to be the two eigenvalues for $C_i$,
we work with the diagram
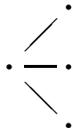
\[
\begin{picture}(16,31)
\put(0,15){\circle*{1}}
\put(15,0){\circle*{1}}
\put(15,15){\circle*{1}}
\put(15,30){\circle*{1}}
\put(12,3){\line(-1,1){8}}
\put(12,15){\line(-1,0){8}}
\put(12,27){\line(-1,-1){8}}
\end{picture}
\]
Since none of the $C_i$ consists of a multiple of the
identity matrix, $\alpha$ is the real root
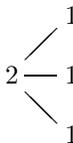
\[
\begin{picture}(18,40)
\put(0,15){2}
\put(15,0){1}
\put(15,15){1}
\put(15,30){1}
\put(13,5){\line(-1,1){8}}
\put(13,17){\line(-1,0){8}}
\put(13,29){\line(-1,-1){8}}
\end{picture}
\]
(where we display the components of $\alpha$ at the appropriate vertices in the diagram).
The sum of all six eigenvalues being zero corresponds to the obvious
requirement that $\tr(A_1)+\tr(A_2)+\tr(A_3)=0$ for $A_i\in C_i$,
and also to the condition $\lambda\cdot\alpha=0$.
Now, for example, the decomposition
\[
\begin{picture}(108,40)
\put(0,15){2}
\put(15,0){1}
\put(15,15){1}
\put(15,30){1}
\put(13,5){\line(-1,1){8}}
\put(13,17){\line(-1,0){8}}
\put(13,29){\line(-1,-1){8}}
\put(30,15){$=$}
\put(45,15){1}
\put(60,0){1}
\put(60,15){0}
\put(60,30){0}
\put(58,5){\line(-1,1){8}}
\put(58,17){\line(-1,0){8}}
\put(58,29){\line(-1,-1){8}}
\put(75,15){$+$}
\put(90,15){1}
\put(105,0){0}
\put(105,15){1}
\put(105,30){1}
\put(103,5){\line(-1,1){8}}
\put(103,17){\line(-1,0){8}}
\put(103,29){\line(-1,-1){8}}
\end{picture}
\]
corresponds to the condition $\xi_{1,1}+\xi_{2,1}+\xi_{3,2}\neq 0$,
or equivalently $\xi_{1,2}+\xi_{2,2}+\xi_{3,1}\neq 0$.

As another example, consider 3 conjugacy classes $C_1$, $C_2$, $C_3$ of diagonalizable
$3\times 3$ matrices, where $C_1$ and $C_2$ have distinct eigenvalues $\xi_{i,1},\xi_{i,2},\xi_{i,3}$ ($i=1,2$),
and $C_3$ has eigenvalues $\xi_{3,1},\xi_{3,2}$ of multiplicities 2, 1
respectively. Thus we deal with the diagram and vector $\alpha$
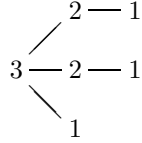
\[
\begin{picture}(33,40)
\put(0,15){3}
\put(15,0){1}
\put(15,15){2}
\put(15,30){2}
\put(30,15){1}
\put(30,30){1}
\put(13,5){\line(-1,1){8}}
\put(13,17){\line(-1,0){8}}
\put(13,29){\line(-1,-1){8}}
\put(28,17){\line(-1,0){8}}
\put(28,32){\line(-1,0){8}}
\end{picture}
\]
Since $\alpha$ is a real root, there is a rigid solution
for general $\lambda$ satisfying $\lambda\cdot\alpha=0$,
or equivalently, for general $\xi_{i,j}$ satisfying
\begin{equation}
\label{eqall}
\xi_{1,1}+\xi_{1,2}+\xi_{1,3}+\xi_{2,1}+\xi_{2,2}+\xi_{2,3}+2\xi_{3,1}+\xi_{3,2}=0.
\end{equation}
There is no rigid solution when
\begin{equation}
\label{eqnosol}
\xi_{1,1}+\xi_{2,1}+\xi_{3,1}=0,
\end{equation}
because of the decomposition
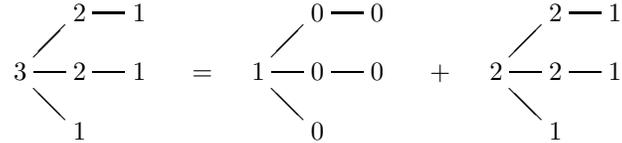
\[
\begin{picture}(153,40)
\put(0,15){3}
\put(15,0){1}
\put(15,15){2}
\put(15,30){2}
\put(30,15){1}
\put(30,30){1}
\put(13,5){\line(-1,1){8}}
\put(13,17){\line(-1,0){8}}
\put(13,29){\line(-1,-1){8}}
\put(28,17){\line(-1,0){8}}
\put(28,32){\line(-1,0){8}}
\put(45,15){$=$}
\put(60,15){1}
\put(75,0){0}
\put(75,15){0}
\put(75,30){0}
\put(90,15){0}
\put(90,30){0}
\put(73,5){\line(-1,1){8}}
\put(73,17){\line(-1,0){8}}
\put(73,29){\line(-1,-1){8}}
\put(88,17){\line(-1,0){8}}
\put(88,32){\line(-1,0){8}}
\put(105,15){$+$}
\put(120,15){2}
\put(135,0){1}
\put(135,15){2}
\put(135,30){2}
\put(150,15){1}
\put(150,30){1}
\put(133,5){\line(-1,1){8}}
\put(133,17){\line(-1,0){8}}
\put(133,29){\line(-1,-1){8}}
\put(148,17){\line(-1,0){8}}
\put(148,32){\line(-1,0){8}}
\end{picture}
\]
where both terms in the sum are roots.
On the other hand, there can be rigid solutions when
\begin{equation}
\label{eqissol}
\xi_{1,1}+\xi_{2,1}+\xi_{3,2}=0.
\end{equation}
Namely, decompositions of $\alpha$ of the form
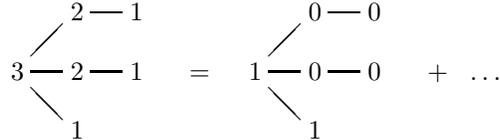
\[
\begin{picture}(115,40)
\put(0,15){3}
\put(15,0){1}
\put(15,15){2}
\put(15,30){2}
\put(30,15){1}
\put(30,30){1}
\put(13,5){\line(-1,1){8}}
\put(13,17){\line(-1,0){8}}
\put(13,29){\line(-1,-1){8}}
\put(28,17){\line(-1,0){8}}
\put(28,32){\line(-1,0){8}}
\put(45,15){$=$}
\put(60,15){1}
\put(75,0){1}
\put(75,15){0}
\put(75,30){0}
\put(90,15){0}
\put(90,30){0}
\put(73,5){\line(-1,1){8}}
\put(73,17){\line(-1,0){8}}
\put(73,29){\line(-1,-1){8}}
\put(88,17){\line(-1,0){8}}
\put(88,32){\line(-1,0){8}}
\put(105,15){$+\ \ \dots$}
\end{picture}
\]
could perhaps fail condition (iii) because the dot product of $\lambda$
with the indicated summand is zero, but the difference
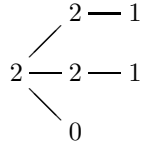
\[
\begin{picture}(33,40)
\put(0,15){2}
\put(15,0){0}
\put(15,15){2}
\put(15,30){2}
\put(30,15){1}
\put(30,30){1}
\put(13,5){\line(-1,1){8}}
\put(13,17){\line(-1,0){8}}
\put(13,29){\line(-1,-1){8}}
\put(28,17){\line(-1,0){8}}
\put(28,32){\line(-1,0){8}}
\end{picture}
\]
is not a root, so this decomposition of $\alpha$ must have at least 3 summands,
and the dot product of $\lambda$ with any of the later summands will be nonzero
for general $\xi_{i,j}$ satisfying (\ref{eqall}) and (\ref{eqissol}).
(I am grateful to a referee of an earlier version of the paper for suggesting this example.)

Using the methods developed in \cite{CBce,CBmm,CBH},
we now give short proof that conditions (i), (ii) and (iii) characterize
the existence of a unique simple representation of $\Pi^\lambda$,
avoiding Theorem \ref{t:quoted}, whose proof in \cite{CBmm} is
quite complicated. Assuming either that there is a unique simple
of dimension vector $\alpha$, or that $\alpha$ satisfies (i)-(iii) we
deduce the other.

Since the relations for $\Pi^\lambda$ are given by a moment map,
a standard symplectic geometry argument shows that the simple representations
of $\Pi^\lambda$, if any, depend on $2p(\alpha)$ parameters.
(For example, using the notation of \cite[Theorem 10.3]{CBce}, the simple representations
form an open subset $\mathcal{S}$ of the
vector space $\operatorname{Rep}(K\overline{Q},\alpha)$ of representations
of $\overline{Q}$ of dimension vector $\alpha$.
Representations of $\Pi^\lambda$ correspond to a fibre of
$\mu:\operatorname{Rep}(K\overline{Q},\alpha)\to \operatorname{End}(\alpha)_0$,
whose restriction to $\mathcal{S}$ is smooth,
and isomorphism classes correspond to orbits of the
group $\operatorname{GL}(\alpha)$ which acts on $\mathcal{S}$ with one-dimensional stabilizers.
Thus the number of parameters is
$\dim \mathcal{S} - \dim \operatorname{End}(\alpha)_0 - (\dim \operatorname{GL}(\alpha) - 1) = 2p(\alpha)$,
as claimed.)

Thus either condition implies that $p(\alpha)=0$, so $\alpha^T \Ca \alpha > 0$,
and hence $\epsilon_v^T \Ca \alpha > 0$ for some vertex $v$. Thus the reflection $s_v(\alpha)$
is strictly smaller than $\alpha$.

Now if $\lambda_v\neq 0$ then the reflection functors of \cite{CBH} show that
simple representations of $\Pi^\lambda$ of dimension $\alpha$ correspond to simple
representations of $\Pi^{\lambda'}$ of dimension $s_v(\alpha)$ for suitable $\lambda'$. Moreover
conditions (i)-(iii) for $\alpha$ correspond to the equivalent conditions
for $s_v(\alpha)$ by \cite[Lemma 5.2]{CBmm}. Thus the claim follows by an induction.

On the other hand, if $\lambda_v=0$ then there is a unique simple representation
if and only if $\alpha=\epsilon_v$ by \cite[Lemma 7.2]{CBmm}, and conditions (i)-(iii)
hold if and only if $\alpha=\epsilon_v$, for otherwise, either (i) or (ii) fails, or
the decomposition $\alpha = s_v(\alpha)+\epsilon_v+\dots+\epsilon_v$ contradicts (iii).

Finally, we remark that Dettweiler and Reiter's additive version \cite[Appendix]{DR}
of Katz's algorithm
corresponds to the the reflection functor at the vertex 0 for the generalized
Cartan matrix introduced in the introduction. (The reflection functors
at other vertices have the same effect as re-ordering the elements $\xi_{i,j}$.)
Either approach can in principle be used to construct the solutions in the rigid
cases. See \cite{G} for an alternative explicit construction of solutions in
some special cases.

\section{The nilpotent case}
\label{s:nilpotent}
The case where the $C_i$ are nilpotent conjugacy classes
has been considered by Kostov \cite{Kostovmon,Kostovaspects}.
Here we may take all $\xi_{i,j}=0$, so that $\lambda=0$.
Now \cite[Theorem 8.1]{CBmm} characterizes the elements of
$\Sigma_0$ as the coordinate vectors, and
the elements of the fundamental region, except those of
three special types (I), (II) and (III). For the generalized Cartan matrix constructed in
the introduction, the fundamental region consists of the nonzero vectors $\alpha\in\N^I$
with connected support, $2\alpha_{[i,j]} \le \alpha_{[i,j-1]}+\alpha_{[i,j+1]}$
for $1\le i\le k$ and $1\le j<d_i-1$,
$2\alpha_{[i,d_i-1]} \le \alpha_{[i,d_i-2]}$ for $1\le i\le k$, and
\[
2\alpha_0 \le \sum_{i=1}^k \alpha_{[i,1]}.
\]
For $\alpha$ constructed as in Theorem \ref{t:main},
the first conditions follow automatically from (\ref{e:rcond}),
and the last condition becomes
\[
2n \le \sum_{i=1}^k r_{i,1}.
\]
Deleting any vertices $v$ at which $\alpha_v=0$, the special types are:
\begin{itemize}
\item[(I)]The diagram is an extended Dynkin diagram $\Delta$, and $\alpha$ is a proper
multiple of the minimal positive imaginary root $\delta$.

\item[(II)]The diagram is obtained from an extended Dynkin diagram $\Delta$ by adding
a new vertex $w$ connected to an extending vertex of $\Delta$. The
restriction of $\alpha$ to $\Delta$ is a proper multiple of $\delta$, and
$\alpha_w=1$.

\item[(III)]Does not occur.
\end{itemize}

As an example, suppose that $C_1$, $C_2$ and $C_3$ are all the conjugacy class of nilpotent $12\times 12$
matrices consisting of four $3\times 3$ Jordan blocks. In this case the diagram and
vector $\alpha$ are
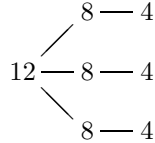
\[
\begin{picture}(33,40)
\put(-3,15){12}
\put(15,0){8}
\put(15,15){8}
\put(15,30){8}
\put(30,0){4}
\put(30,15){4}
\put(30,30){4}
\put(13,5){\line(-1,1){8}}
\put(13,17){\line(-1,0){8}}
\put(13,29){\line(-1,-1){8}}
\put(28,2){\line(-1,0){8}}
\put(28,17){\line(-1,0){8}}
\put(28,32){\line(-1,0){8}}
\end{picture}
\]
This is of type (I), with the diagram being extended Dynkin of type $\tilde{E}_6$, so there is no irreducible
solution. Changing $C_3$ to involve Jordan blocks of sizes $4,3,3,2$, we obtain
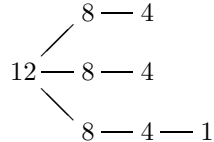
\[
\begin{picture}(48,40)
\put(-3,15){12}
\put(15,0){8}
\put(15,15){8}
\put(15,30){8}
\put(30,0){4}
\put(30,15){4}
\put(30,30){4}
\put(45,0){1}
\put(13,5){\line(-1,1){8}}
\put(13,17){\line(-1,0){8}}
\put(13,29){\line(-1,-1){8}}
\put(28,2){\line(-1,0){8}}
\put(28,17){\line(-1,0){8}}
\put(28,32){\line(-1,0){8}}
\put(43,2){\line(-1,0){8}}
\end{picture}
\]
which is of type (II), so again there is no irreducible solution. Changing $C_3$ again, to involve
Jordan blocks of sizes $4,4,2,2$, we obtain
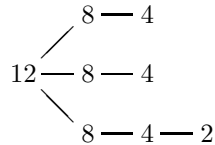
\[
\begin{picture}(48,40)
\put(-3,15){12}
\put(15,0){8}
\put(15,15){8}
\put(15,30){8}
\put(30,0){4}
\put(30,15){4}
\put(30,30){4}
\put(45,0){2}
\put(13,5){\line(-1,1){8}}
\put(13,17){\line(-1,0){8}}
\put(13,29){\line(-1,-1){8}}
\put(28,2){\line(-1,0){8}}
\put(28,17){\line(-1,0){8}}
\put(28,32){\line(-1,0){8}}
\put(43,2){\line(-1,0){8}}
\end{picture}
\]
which is in the fundamental region (since $2\times 12 \le 8+8+8$), and is not a special type,
so there is an irreducible solution in this case.

Observe that type (I) corresponds to the `special' cases,
and type (II) to the `almost special' cases (a1), (b1), (c2) and (d3) of \cite[\S 6]{Kostovaspects}.
Thus one implication in our characterization has already been obtained by
Kostov \cite[Theorem 34]{Kostovaspects}.

\section{The generic eigenvalues case}
\label{s:generic}
Let $\nu_{i,1},\dots,\nu_{i,r_i}$ be the distinct
eigenvalues of $C_i$, and let $m_{i,\ell}$ be the
algebraic multiplicity of $\nu_{i,\ell}$.
We assume that
\begin{equation}
\label{e:geneq}
\sum_{i=1}^k \sum_{\ell=1}^{r_i} m_{i,\ell} \nu_{i,\ell}= 0,
\end{equation}
which is an obvious necessary condition for equation~(\ref{eq2}) to have
a solution. Following Kostov~\cite{KostovCR,Kostovnext}, we say that the $C_i$ have
\emph{generic eigenvalues} if, for any integers $0 \le m'_{i,\ell} \le m_{i,\ell}$
such that $\sum_{\ell=1}^{r_i} m'_{i,\ell}$ is independent of $i$, an equality
of the form
\[
\sum_{i=1}^k \sum_{\ell=1}^{r_i} m'_{i,\ell} \nu_{i,\ell}= 0
\]
implies that $m'_{i,\ell}=0$ for all $i,\ell$, or $m'_{i,\ell}=m_{i,\ell}$ for all $i,\ell$.
Clearly, if the $C_i$ have generic eigenvalues, then
the $m_{i,\ell}$ have no common divisor.

In \cite{KostovCR,Kostovnext}, Kostov determines those $C_i$ with generic
eigenvalues for which equation (\ref{eq2}) has an irreducible solution.
His answer can be neatly reformulated in terms of roots,
and we show that it also follows from Theorem~\ref{t:main}.
Choose elements $\xi_{i,j}$ as in the introduction, and
let $\lambda$ and $\alpha$ be constructed as in Theorem~\ref{t:main}.
Assuming that the $C_i$ have generic eigenvalues, we show that $\alpha\in\Sigma_\lambda$
(so that equation~(\ref{eq2}) has an irreducible solution) if and only
if $\alpha$ is a root.

First note that $\alpha$ is indivisible, that is, its components have no common
divisor, for otherwise the $m_{i,\ell}$ will have the same common divisor.
Also, equation (\ref{e:geneq}) implies that $\lambda\cdot\alpha=0$.

By definition, if $\alpha\in\Sigma_\lambda$, then $\alpha$ is a root.
For the converse, suppose that $\alpha$ is a root, and that
$\alpha = \beta^{(1)}+\beta^{(2)}+\dots$
is a nontrivial decomposition of $\alpha$ as a sum of elements of $R_\lambda^+$.
We show that inequality (\ref{e:pmax}) holds.

Let $\beta^{(s)}$ be a summand occurring in this decomposition. For given $i$
and $\ell$, define
\[
m^{(s)}_{i,\ell} = \sum_{\substack{j=1 \\ \xi_{i,j}=\nu_{i,\ell}}}^{d_i}
(\beta^{(s)}_{[i,j-1]} - \beta^{(s)}_{[i,j]}),
\]
where, for convenience, we define $\beta^{(s)}_{[i,d_i]}=0$,
and $[i,0]$ denotes the vertex 0. Clearly
\[
\sum_{\ell=1}^{r_i} m^{(s)}_{i,\ell} = \sum_{j=1}^{d_i} (\beta^{(s)}_{[i,j-1]} - \beta^{(s)}_{[i,j]}) = \beta^{(s)}_0,
\]
which is independent of $i$.

We claim that $m^{(s)}_{i,\ell} \ge 0$ for all $i,\ell$. To prove it, we divide into
two cases.
If $\beta^{(s)}_0\neq 0$, then the fact that $\beta^{(s)}$ is a root implies that
$\beta^{(s)}_{[i,j-1]} \ge \beta^{(s)}_{[i,j]}$ for all $i,j$, which immediately
gives the claim.
On the other hand, if $\beta^{(s)}_0=0$, then since any root has connected support,
$\beta^{(s)}$ is supported on one arm of the diagram, say the $i$-th. Moreover,
using the classification of roots for type $A_n$, it follows that
there are $1\le p\le q\le d_i-1$ with
\begin{equation}
\label{e:armroot}
\beta^{(s)}_{[i,p]} = \beta^{(s)}_{[i,p+1]} = \dots = \beta^{(s)}_{[i,q]} = 1,
\end{equation}
and all other components of $\beta^{(s)}$ are zero. Now
\[
0 = \lambda\cdot\beta^{(s)} = \lambda_{[i,p]} + \lambda_{[i,p+1]} + \dots + \lambda_{[i,q]} = \xi_{i,p}-\xi_{i,q+1}.
\]
so that $\xi_{i,p}=\xi_{i,q+1}$. From this it follows that $m^{(s)}_{i,\ell}=0$ for all $\ell$, and clearly also
$m^{(s)}_{i,\ell}=0$ for all other arms $i$.

Now, varying $s$, we have
\[
m^{(1)}_{i,\ell} + m^{(2)}_{i,\ell} + \dots =
\sum_{\substack{j=1 \\ \xi_{i,j}=\nu_{i,\ell}}}^{d_i} (\alpha_{[i,j-1]} - \alpha_{[i,j]}) = m_{i,\ell}
\]
(where, again, we define $\alpha_{[i,d_i]}=0$).
Thus $0 \le m^{(s)}_{i,\ell} \le m_{i,\ell}$.
Moreover,
\begin{align*}
0 &= \lambda\cdot\beta^{(s)} =
\sum_{i,j} (\xi_{i,]}-\xi_{i,j+1}) \beta^{(s)}_{[i,j]} - (\sum_i \xi_{i,1}) \beta^{(s)}_0
\\
&= \sum_{i,j} \xi_{i,j} (\beta^{(s)}_{[i,j]} - \beta^{(s)}_{[i,j-1]})
= - \sum_{i,\ell} \nu_{i,\ell}\, m^{(s)}_{i,\ell}.
\end{align*}
This contradicts the genericity assumption unless $m^{(s)}_{i,\ell}=m_{i,\ell}$ for all $i,\ell$,
or $m^{(s)}_{i,\ell}=0$ for all $i,\ell$.
Thus, exactly one of the summands $\beta^{(s)}$ has $\beta^{(s)}_0\neq 0$.
By renumbering we may suppose it is the one with $s=1$.

We show that $\alpha^T \Ca \beta^{(s)}\le 0$ for $s\neq 1$, where
$\Ca$ is the generalized Cartan matrix.
Namely, using that $\beta^{(s)}$ is of the form (\ref{e:armroot}), we have
\[
\alpha^T \Ca \beta^{(s)} = (\alpha_{[i,q]} - \alpha_{[i,q+1]}) - (\alpha_{[i,p-1]}-\alpha_{[i,p]}),
\]
so the assertion follows from equation (\ref{e:rcond}).

Now let $\gamma = \beta^{(2)}+\beta^{(3)}+\dots = \alpha-\beta^{(1)} \neq 0$.
We have $\gamma^T \Ca \gamma > 0$ since the support of $\gamma$ is contained
in the diagram obtained by deleting the vertex 0, which is a union of Dynkin
diagrams of type $A_n$, so has positive definite quadratic form.
Since $\Ca$ is symmetric we have
\begin{align*}
p(\beta^{(1)}) &= 1- \frac12 (\beta^{(1)})^T \Ca \beta^{(1)}
= 1- \frac12 \alpha^T \Ca \alpha - \frac12 \gamma^T \Ca \gamma + \alpha^T \Ca \gamma \\
&= p(\alpha)- \frac12 \gamma^T \Ca \gamma + \alpha^T \Ca \beta^{(2)}+ \alpha^T \Ca \beta^{(3)}+\dots
< p(\alpha).
\end{align*}
Since also $p(\beta^{(s)}) = 0$ for $s\neq 1$, inequality (\ref{e:pmax}) follows.
Thus $\alpha\in\Sigma_\lambda$.


\begin{thebibliography}{99}
\bibitem{AB}
D. V. Anosov and A. A. Bolibruch,
\emph{The Riemann-Hilbert problem},
Friedr.\ Vieweg \& Sohn, Braunschweig, 1994.

\bibitem{CBce}
W. Crawley-Boevey,
\emph{Preprojective algebras, differential operators and a Conze embedding for deformations of Kleinian singularities},
Comment. Math.\ Helv., \textbf{74} (1999), 548--574.

\bibitem{CBmm}
W. Crawley-Boevey,
\emph{Geometry of the moment map for representations of quivers},
Compositio Math., \textbf{126} (2001), 257--293.

\bibitem{CBH}
W. Crawley-Boevey and M. P. Holland,
\emph{Noncommutative deformations of Kleinian singularities},
Duke Math.\ J., \textbf{92} (1998), 605--635.

\bibitem{DR}
M. Dettweiler and S. Reiter,
\emph{An algorithm of Katz and its application to the inverse Galois problem},
J.\ Symbolic Computation \textbf{30} (2000), 761--798.

\bibitem{G}
O. A. Gleizer,
\emph{Explicit Solutions to the Additive Deligne-Simpson Problem and Their Applications},
preprint math.LA/0105184.

\bibitem{Katz}
N. M. Katz,
\emph{Rigid local systems},
Princeton University Press, Princeton, NJ, 1996.

\bibitem{Kostovmon}
V. P. Kostov,
\emph{On the existence of monodromy groups of Fuchsian systems
on Riemann's sphere with unipotent generators},
J. Dynam.\ Control Systems, \textbf{2} (1996), 125--155.

\bibitem{KostovCR}
V. P. Kostov,
\emph{On the Deligne-Simpson problem},
C. R. Acad.\ Sci.\ Paris S\'er.\ I Math.\ \textbf{329} (1999), 657--662.

\bibitem{Kostovaspects}
V. P. Kostov,
\emph{On some aspects of the Deligne-Simpson problem},
preprint math.AG/0005016.

\bibitem{Kostovnext}
V. P. Kostov,
\emph{On the Deligne-Simpson problem},
preprint math.AG/0011013.

\bibitem{Simp}
C. T. Simpson,
\emph{Products of Matrices},
Differential geometry, global analysis, and topology (Halifax, NS, 1990),
Canadian Math.\ Soc.\ Conf.\ Proc.\ \textbf{12} (1992),
Amer.\ Math.\ Soc., Providence, RI, 1991,
pp.\ 157--185.

\end{thebibliography}
\end{document}